\documentclass[letterpaper, 10 pt, conference]{ieeeconf} 
\IEEEoverridecommandlockouts 
\usepackage{textcomp}
\usepackage{xr-hyper}
\usepackage[utf8]{inputenc}
\usepackage{graphicx} 
\usepackage{cite}
\usepackage{times} 
\usepackage{amsmath} 
\usepackage{amssymb}  
\usepackage{lipsum}
\usepackage{color}

\let\labelindent\relax
\usepackage{amsthm}
\usepackage{enumitem}
\usepackage{comment}
\usepackage{multirow}
\usepackage{algorithm}
\makeatletter
\renewcommand{\ALG@name}{Step}
\makeatother
\usepackage{algpseudocode}
\usepackage{authblk}
\usepackage[normalem]{ulem}
\usepackage{cases}
 \usepackage{mathtools}
 \usepackage{nccmath}
\algnewcommand\algorithmicforeach{\textbf{for each}}
\algdef{S}[FOR]{ForEach}[1]{\algorithmicforeach\ #1\ \algorithmicdo}
\newcommand{\stkout}[1]{\ifmmode\text{\sout{\ensuremath{#1}}}\else\sout{#1}\fi}
\definecolor{mypink}{RGB}{255,56,152}

\definecolor{mygreen}{RGB}{23,198,38}

\newcommand{\Spanpos}{\mathrm{Span}^+}
\newcommand{\V}{\mathcal{V}}

\newcommand{\myL}{\mathcal{L}}
\newcommand{\E}{\mathcal{E}}
\newcommand{\D}{\mathcal{D}} 
\newcommand{\Real}{\mathbb{R}}
\newcommand{\Span}{\mathrm{Span}}
\newcommand{\proj}{\mathrm{proj}}
\newcommand{\G}{\mathcal{G}} 
\newcommand{\C}{\mathcal{C}} 
\newcommand{\K}{\mathcal{K}} 
\newcommand{\T}{\mathrm{T}} 
\newcommand{\W}{\mathcal{W}} 
\newcommand{\Chi}{\mathcal{X}}
\newcommand{\Q}{\mathcal{Q}} 
\newcommand{\e}{\mathfrak{e}}
\DeclareMathOperator*{\argmax}{arg\,max} 
\DeclareMathSymbol{\shortminus}{\mathbin}{AMSa}{"39}
\usepackage{graphicx}
\makeatletter
\DeclareMathOperator{\sgn}{sgn}
\newcommand{\shorteq}{%
  \settowidth{\@tempdima}{-}
  \resizebox{\@tempdima}{\height}{=}%
}
\makeatother
\newtheorem{mydef}{Definition}
\newtheorem{mylem}{Lemma}
\newtheorem{myrem}{Remark}
\newtheorem{mythm}{Theorem}
\newtheorem{mycor}{Corollary}

\newtheorem{assumption}{Assumption}

\begin{document}
\title{Partial controllability of network dynamical systems with unilateral inputs}
\author{Camilla Ancona, Francesco Lo Iudice, Antonio Coppola,  Pietro De Lellis and Franco Garofalo
\thanks{All authors are with the Department of Information Technology and Electrical Engineering, University of Naples Federico II, Naples, Italy.}
\thanks{This work was supported by the program “STAR 2018” of the University of Naples Federico II and Compagnia di San Paolo, Istituto Banco di Napoli - Fondazione, project ACROSS and by the Italian Ministry of University and Research (2020–2023) under the Research Project PRIN 2017 “Advanced Network Control of Future Smart Grids”.\\
© 2022 IEEE.  Personal use of this material is permitted.  Permission from IEEE must be obtained for all other uses, in any current or future media, including reprinting/republishing this material for advertising or promotional purposes, creating new collective works, for resale or redistribution to servers or lists, or reuse of any copyrighted component of this work in other works.}
}%

\maketitle
\thispagestyle{empty}
\pagestyle{empty}
\begin{abstract}
    Our ability to control network dynamical systems is often hindered by constraints on the number and nature of the available control actions, which make controlling the whole network unfeasible. In this manuscript, we focus on the case where unilateral inputs are exerted on a subset of the network nodes. Leveraging the observation that, different from the case of subsystems, unilateral node reachability and controllability are equivalent, we provide conditions for a given node subset to be unilaterally controllable. The theoretical findings are then employed to develop a computationally efficient heuristic to select the nodes where the unilateral inputs should be injected.
\end{abstract}
\section{Introduction} \label{intro}
Modeling, analyzing, and controlling network dynamical systems is of interest for applications as diverse as formation control \cite{delellis2015decentralised}, multicellular control in biology \cite{fiore2020multicellular,iglesias2010control}, power systems \cite{hommelberg2007distributed}, and financial market dynamics \cite{caldarelli2004emergence,de2018overconfident}. In the last decades, control engineers have focused on designing distributed protocols capable of inducing the emergence of collective behaviors, such as e.g. consensus and synchronization \cite{yu2010some, wu1995synchronization,song2010second,li2009consensus}. More recently, the ambition to systematically tackle more general network control problems has brought to specify the classical concept of controllability to the case of network systems \cite{liu2011controllability}. It has been pointed out that, when studying network dynamical systems, (i) controllability should be conferred through an appropriate choice of the nodes where the inputs should be injected, thereby several input selection algorithms have been developed \cite{gao2014target,wang2017physical}, (ii) existing controllability tests may be misleading since, when the number of inputs is much smaller than the number of network nodes, controlling a network can be energetically unfeasible \cite{pasqualetti2014controllability,lindmark2018minimum,de2021inherent}, and (iii) achieving complete network controllability can turn out being unnecessary or unfeasible, and thus one should rather focus on controlling selected subnetworks \cite{delellis2018partial}. 

When studying controllability, a crucial difference exists between large scale dynamical systems and network systems, whereby subsytem controllability differs from controllability of a node subset for a subtle, yet critical, aspect \cite{LiGaSo:19}. Although in general the choice of the reference frame for the state variables of the controllable subsystem is irrelevant, this is not true for node subsets, as we need to preserve the association between nodes and network state variables. This in turn has several relevant consequences, such as the fact that, whereas the controllable subsystem is unique, there can be multiple controllable node subsets. 

A further challenge, typically overlooked in the literature on network control, is that in real-world systems the input signals are constrained. A relevant case in applications is when the inputs are constrained to be unilateral, that is, the sign of each signal cannot change over time\cite{goodwine1996controllability,LinAlt:17}. Practical examples where such constraints arise include the optimal power flow problem in power grids, where nodes are either loads or generators \cite{frank2012optimal,frank2012optimal2}, the inhibition or activation or genes in biological networks \cite{mcdonald2008activating}, control of wire-driven parallel robots \cite{merlet2013wire,alamdari2018stiffness}, or marketing campaigns where comparative advertising is forbidden \cite{miskolczi2004definition}, see \cite[Table 1]{LinAlt:17} for further practical instances of unilateral control. The literature on controllability under constrained inputs can be traced back to the Seventies \cite{saperstone1971controllability,Brammer:72}, but only recently the problem has been tackled for network dynamical systems. Specifically, Lindmark and Altafini have derived conditions for finding the minimal set of inputs that render the whole network controllable  \cite{LinAlt:17}.

To the best of our knowledge, none of the existing work tackled the controllability problem of a subset of the network nodes with unilateral inputs, and this is the gap we aim at filling in this manuscript. Different from the case of unconstrained inputs, when the inputs are constrained to be unilateral, the controllability of a network cannot be studied through structural approaches, nor through the controllability gramian. This in turn restricts the theoretical tools available to design optimal input placement strategies, making the unilateral case much more challenging than the unconstrained case. We tackle this problem by first characterizing the convex cone containing the unilaterally controllable states of a linear dynamical system, which we show to differ from the one containing the unilaterally reachable states. Then, by means of a suitable projection, we translate these results for network systems, and obtain the conditions for unilateral reachability and controllability of a node subset that, different from the general case of systems, we observe to be equivalent. The theoretical findings are then used to develop a greedy heuristic to decide where to inject the unilateral inputs, which provides a suboptimal solution to the problem of maximizing the number of controllable nodes.

\section{Preliminaries} \label{prel}

Given a set $\mathcal H$, we denote by $|\mathcal H|$ its cardinality, and given a vector space, we denote by $0$ its origin. Given a real vector space $\Real^n$, we denote by $\Real^n_{\ge0}$ ($\Real^n_{\le 0}$) the set of vectors in $\Real^n$ with nonnegative (nonpositive) entries.
Let $\D$ be a set of $|\D|=k$ vectors $d_1,\dots,d_k$ in $\Real^n$, $\Span(\D)$ is the set of all linear combinations of the vectors in $\D$. The \emph{positive span} $\Spanpos(\D)$ of ${\D}$ is the set of all linear combinations with nonnegative coefficients, that is, ${\Spanpos(\D)} = \{\sum_{i=1}^{k}\alpha_id_i \;:\;\alpha_i \in \Real_{\geq 0}\}$, which constitutes a polyhedral convex cone \cite{schrijver1998theory}.
If $\D$ is a singleton then $\Spanpos(\D)$ is called a ray. All the rays and the singleton $\mathcal{O}=\{ 0\}$ are degenerate cones.
Given a convex polyhedral cone $\C$, we define its dimension $|\C|$ as the number of vectors required to generate it. 
The lineality space of a convex cone $\C$ is defined as the largest subspace $\Chi^l:= \C \cap -\C $ contained in $\C$, whose dimension is the lineality of $\C$ \cite{onconvexcones}.

Next, let us denote by $\e_i$ the $i$-th versor in $\Real^n$. Given an index set $\K$, we define $\Chi_\K$ as the subspace linearly spanned by $\cup_{i \in \K} \{\e_i\}$. Furthermore, given a vector $d\in\mathbb{R}^n$, we denote by ${\proj}_{\Chi_\K}(d) = \sum_{{i \in \K}}\big(d^{\T}\e_i\big)\e_i$ the orthogonal projection of $d$ along $\Chi_\K$. Given a complex vector $c\in\mathbb C^n$, we denote by $\Re({c})$ and $\Im({c})$ its real and imaginary parts, respectively. The operators $\lor$ and $\land$ denote the logical disjunction and conjunction, respectively, whereas the symbol $\setminus$ denotes a set difference. Finally, the big-O notation $O(\cdot)$ 
describes the order of magnitude of the algorithm execution time with respect to the number of steps required to complete it.

\section{Problem formulation}\label{sec:pr_form}
Let us consider a linear dynamical network on a graph $\mathcal G =\{ \mathcal V, \mathcal E\}$, where $\V$ and $\E \subseteq (\V \times \V)$ are the sets of its nodes and edges, respectively. Defining the network state $x=[x_1,\ldots,x_n]^\T$, with $x_i\in\mathbb{R}$ being the state of the $i$-th node, the network dynamics are given by
\begin{equation} \label{eq:net_eq}
\dot{x}(t)=Ax(t)+Bu(t),
\end{equation}
where $A\in\mathbb{R}^{n\times n}$ is the adjacency matrix of $\mathcal G$, whose $ij$-th entry $a_{ij}\neq 0$ if $(i,j)\in\mathcal{E}$ and $a_{ij}=0$ otherwise. Matrix $A$ encapsulates both the individual dynamics and the interaction between the network nodes, which are encoded by the diagonal and off-diagonal elements of $A$, respectively. Matrix $B\in\mathbb{R}^{n\times m}$ is the input matrix, whose $ij$-th element modulates the effect the input $u_j$ has on the dynamics of node $i$. Here, we consider the case of unilateral inputs.
\begin{mydef}
The input $u(t)$ to network \eqref{eq:net_eq} is called unilateral if $(u(t)\in \Real_{\ge 0}^m ,\forall t) \lor (u(t)\in \Real_{\le 0}^m ,\forall t)$.
\end{mydef}
In what follows, without loss of generality, we consider only nonnegative inputs for the sake of clarity, whereby their sign will be determined by the sign of the entries of $B$. More formally, we make the following assumption.
\begin{assumption}\label{assump:unilateral}
Inputs are nonnegative, that is $u(t)\in\Real^m_{\ge 0}$ for all $t$ and each column of matrix $B$ belongs to the set $\mathcal{B} = \{\e_i, -\e_i,i=1,\ldots,n\}$ \cite{LinAlt:17}. \end{assumption}
Here, we focus on the case in which unilateral controllability of the whole network is not feasible, whereby the conditions given in \cite{LinAlt:17} do not hold. The problem then arises of selecting the input so that the state of a subset of the network nodes can be steered towards any desired value. Before stating this problem, we need to define unilateral controllability of a node subset $\V_s$, whose associated state $x_s$ is the vector stacking the states of all nodes in $\V_s$.  

\begin{mydef}\label{def:uni_reach}
A node subset $\V_s\subseteq \V$ is unilaterally reachable if the state of its nodes $x_s$ can be steered from $0$ to any target value in finite time through an appropriate selection of the unilateral input $u(t)$.
\end{mydef}
\begin{mydef}\label{def:uni_cont}
A node subset $\V_s\subseteq \V$ is unilaterally controllable if, for all initial conditions $x_s(0)$, the state of its nodes can be steered towards any target value in finite time through an appropriate selection of the unilateral input $u(t)$. 
\end{mydef}
Given the adjacency matrix $A$, the controllability problem we consider is that of designing the input matrix $B$ fulfilling Assumption \ref{assump:unilateral} that maximizes the cardinality of the set $\V_s=\V_s(B)$ of unilaterally controllable nodes, that is,

\begin{subequations}\label{eq:prob_form}
\begin{align}
& \underset{B\subset \mathcal B}{\max} \quad |\V_s|\\
& \text{subject \ to}\nonumber\\
&\sum_{i,j}|b_{ij}| = m \label{m_inputs}\\
& \V_s \text{ unilaterally controllable}
\end{align}
\end{subequations}
Solving \eqref{eq:prob_form} requires finding the conditions such that, given a set of control inputs, a set of nodes is unilaterally controllable, and then devising an input placement algorithm that finds a unilaterally controllable node subset of maximal dimension.
\section{Unilateral reachability and controllability of a node subset}\label{sec:results}

Let $J$ be the Jordan normal form of matrix $A$ and $\mu$ the number of its blocks. We associate to each Jordan block $J_i$, whose size we denote by $\nu_i$, the corresponding eigenvalue $\lambda_i$, for $i=1,\ldots,\mu$.
Next, let $\myL:=\lbrace l_{i,k}, i=1,\dots,\mu, k = 1,\dots,\nu_i \rbrace$ be a set of $\mu$ chains of unit norm linearly independent (generalized) left eigenvectors of $A$ with the maximal number of elements orthogonal to the columns of matrix $B$. 
Finally, let us denote by $T$ the matrix obtained by juxtaposing row-wise the elements of $\myL$. 
We can now define the set $\gamma_{i,k}(l_{i,k}^{\T}B)$ as
\begin{subnumcases}{\gamma_{i,k} \shorteq}
\label{a}
\thickmuskip=0.5\thickmuskip
\varnothing, \qquad\qquad\qquad\quad \ \ \, \text{if } l_{i,k'}^{\T}B = 0+j0 \ \forall k'\geq k \\
\label{b}
\thickmuskip=0\thickmuskip
\hspace{-1.5mm}\{r_{i,k}, \shortminus r_{i,k}\}, \, \, \text{if } \exists  \ k',k'' \geq k : l_{i,k'}^{\T}B \in (\Real^m _{\geq 0 }\setminus\mathcal O) \land \nonumber\\ \qquad\qquad\quad\ \ \,  \hspace{2mm}l_{i,k''}^{\T}B \in \Real^m_{\leq 0}\setminus\mathcal O), \Im{( l_{i,k} )}=0,\\
\label{c}
\thickmuskip=0.5\thickmuskip
\{r_{i,k}\} \,\qquad\quad\ \, \text{if } \exists \ k' \geq k : l_{i,k'}^{\T}B \in ({\Real^m_{\ge 0}\setminus\mathcal O)} \land \nonumber \\ \quad \ \, \nexists \ k'' \geq k : l_{i,k''}^{\T}B \in ({\Real^m_{\le 0}\setminus\mathcal O)}, \Im{(l_{i,k})}=0,\\
\label{d}
\thickmuskip=0.5\thickmuskip
\{\shortminus r_{i,k}\}, \quad\quad\ \, \,  \text{if } \exists \ k' \geq k : l_{i,k'}^{\T}B \in ({\Real^m_{\le 0}\setminus\mathcal O)} \land \nonumber \\ \qquad \nexists\ k'' \geq k | l_{i,k''}^{\T}B \in ({\Real^m_{\ge 0}\setminus\mathcal O)}, \Im{(l_{i,k})}=0,\\
\label{e}
\thickmuskip=0.01\thickmuskip
 \hspace{-1.5mm} \big\lbrace \Re(r_{i,k}),\shortminus \Re(r_{i,k}), \Im{(r_{i,k})}, \shortminus \Im{(r_{i,k})} \big\rbrace, \ \text{otherwise} 
\end{subnumcases} 
for $i= 1, \dots,\mu$. 
Additionally, we denote by $\mathcal{C}_r(B)$ the positive span of the set of all $\gamma_{i,k}$-s, that is, 
\begin{equation}\label{eq:c_r}
\mathcal{C}_r(B):= \Spanpos\big(\bigcup_{i = 1}^{\mu}{\bigcup_{k = 1}^{\nu_i}\{\gamma_{i,k}\}}\big)
\end{equation}
\begin{mythm}\label{thm:uni_con}
If Assumption \ref{assump:unilateral} holds, then \begin{enumerate}[label=(\roman*),align=left, noitemsep]
    \item the cone $\mathcal{C}_r(B)$ is the set of unilaterally reachable states of the pair $(A,B)$;\label{fact1}
    \item the lineality space $\Chi^l$
    of $\mathcal{C}_r(B)$ is the largest unilaterally reachable subspace of the pair $(A,B)$.\label{fact2}
\end{enumerate}
\end{mythm}
\begin{proof}
{\it Statement (i)}:
Let us consider the transformation $z = Tx$. 
As $J = TA T^{-1} $, the dynamics of network \eqref{eq:net_eq} become 
$\dot{z}(t) = J z(t) + T Bu(t).$
By setting $z(0)=0$, we obtain its forced dynamics as $z(t) = \int_0^t \exp({J(t-\tau)})TBu(\tau)d\tau$
or, in scalar form,
\begin{equation}\label{eq:z_ik}
z_{i,k}(t) = \sum_{k'\geq k}\sum_{j = 1}^{m}l^{\T}_{i,k'}b_j\eta_j(u_j) d\tau, 
\end{equation}
for all $i=1,\dots, \mu$, $k =1,\dots, \nu_i$,
where
$$\eta_j(u_j) = \int_{0}^{t}(t-\tau)^{k'-k}   \exp({ \lambda_{i}(t-\tau)}) u_j(\tau) d\tau.$$ 
Since $x=T^{-1}z$, and as the columns of $T^{-1}$ are right generalized eigenvectors of $A$, $z_{i,k}(t)$ represents the dynamics along the right eigenvector $r_{i,k}$ of each Jordan block $J_{i}$, for all $i=1,\ldots,\mu$. Let us now distinguish the case in which $r_{i,k}$ is associated to a real or to a complex eigenvalue\\
 {\it Case (a): $\Im(\lambda_{i})=0$.} From Assumption \ref{assump:unilateral} (i.e., nonnegative inputs), we have that $\eta_j(u_j)\ge 0$ for all $j$. Hence, from \eqref{eq:z_ik} we have that 
\begin{itemize}
    \item if $l_{i,k'}^{\T}B = 0+j0 \ \forall k'\geq k$, then any $\tilde x \in \Span{( r_{i,k}})$ is unreachable;
    \item if there exist $\exists  \ k',k'' \geq k$ such that $  l_{i,k'}^{\T}b_j \geq 0, l_{i,k''}^{\T}b_j \leq 0$, any $\tilde x \in \Span(r_{i,k})$ is unilaterally reachable;
      \item if $\exists \ k' \geq k : l_{i,k'}^{\T} \in ({\Real^m_{\ge 0}\setminus\mathcal O)}$ and $ \nexists \ k'' \geq k : l_{i,k''}^{\T}B \in ({\Real^m_{\le 0}\setminus\mathcal O)}$, then any $\tilde x\in\Spanpos (r_{i,k})$ is unilaterally reachable;
    \item if $\exists \ k' \geq k : l_{i,k'}^{\T}B \in ({\Real^m_{\le 0}\setminus\mathcal O)} $ and $ \nexists\ k'' \geq k : l_{i,k''}^{\T}B \in ({\Real^m_{\ge 0}\setminus\mathcal O)}$, then any $\tilde x\in\Spanpos(\shortminus r_{i,k})$ is unilaterally reachable.
\end{itemize}
{\it Case (b): $\Im(\lambda_i)\ne 0$.} 
As $A$ is a real matrix, each complex eigenvalue will have a complex conjugate. Therefore, the modal dynamics associated to each of the $\nu_i$ pairs of complex conjugate eigenvalues $(\lambda_{i,k},\lambda^*_{i,k})$, $k=1,\dots,\nu_i$ occur along the plane $\widetilde{\Chi}=\Span(\left\lbrace\Re{(r_{i,k})},\Im{(r_{i,k})}\right\rbrace)$ of $\mathbb{R}^n$ and, according to Euler's formula, can be expressed as a sum of sinusoidal functions. Hence, all the states belonging to $\widetilde{\Chi}$ are unilaterally reachable if there exists an index $j$ such that  $\proj_{\widetilde{\Chi}}l^{\T}_{i,k'}b_j \neq 0, \forall k'\geq k$. If, instead, such an index did not exist, then no state would be unilaterally reachable. Finally, considering that (i) if two (or more) states are unilaterally reachable, then any positive combination of the these states is also unilaterally reachable, and (ii) any linear combination involving an unreachable state defines another unreachable state, Statement (i) follows.
\newline {\it Statement (ii).} From Statement (i), no state outside $\C_r(B)$ is reachable. Therefore, Statement (ii) follows.
\end{proof}
\begin{mylem} \label{lem:contr}
Let Assumption 2 hold, then the set of controllable states is $\C_c(B):=\Chi^l(B)\cup \lbrace x\in\Real^n:\shortminus\exp(At)x\in \C_r(B) \rbrace$.
\end{mylem}
\begin{proof}
The thesis follows from the consideration that a point $\bar x $ is controllable if and only if $\shortminus\exp(At)\bar x$ is reachable.
\end{proof}
\subsection{Node subset unilateral reachability}
When studying partial unilateral controllability of network dynamical systems, we need to preserve the association between state variables and network nodes. Therefore, we now provide the following theorems and corollaries characterizing the unilateral reachability and controllability of a node subset.

\begin{mythm} \label{thm:2}
Given network \eqref{eq:net_eq} and a node subset $\V_s \subset \V$, if Assumption \ref{assump:unilateral} holds and $\proj_{\Chi_{\V_s}}(\C_r(B)) = \Chi_{\V_s}$, then $\V_s$ is unilaterally reachable. \label{statement1:thm2}
\end{mythm}
\begin{proof}
From Definition \ref{def:uni_reach}, for a node subset $\V_s$ to be reachable, for all $\bar x_s\in\Real^{|\V_s|}$ and $x(0): x_s(0) = 0$ there must exists a unilateral input $u(t)$ that steers the network towards a state $\bar x$ such that the projection of $\bar x$ on the subspace $\Chi_{\V_s}$ spanned by the versors $\{\e_i\, :\, i \in \V_s\}$, is $\bar x_s$. This is equivalent to the existence of a point $\bar{\bar {x}}\in \C_r(B)$ such that
 \begin{equation}\label{eq:for_proof}
 \proj_{\Chi_{\V_s}}  \bar{\bar {x}} = \bar x_s - \proj_{\Chi_{\V_s}} \exp(At)x(0).    
 \end{equation}
 As from Theorem \ref{thm:uni_con} $\C_r(B)$ is the unilaterally reachable cone, and $\proj_{\Chi_{\V_s}}(\C_r(B)) = \Chi_{\V_s}$ by hypothesis, a point $\bar{\bar{x}}$ fulfilling \eqref{eq:for_proof} exists for all $\bar x_s$ and $x(0): x_s(0) = 0$, and thus the thesis follows.
\end{proof}
Interestingly, we note that the number of unilaterally reachable nodes may be larger than the dimension of the largest unilaterally reachable subspace.
\subsection{Node subset unilateral controllability}
\begin{mythm} \label{thm:3}
Given network \eqref{eq:net_eq}, if Assumption \ref{assump:unilateral} holds, then a node subset $\V_s$ is unilaterally reachable if and only if it is unilaterally controllable. 
\end{mythm}
\begin{proof}
Unilateral controllability of a node subset trivially implies its unilateral reachability, see Definitions \ref{def:uni_reach} and \ref{def:uni_cont}. Hence, let us focus on proving that unilateral reachability of a node implies its unilateral controllability. From Definition \ref{def:uni_cont}, for a node subset $\V_s$ to be unilaterally controllable, for all $\bar x_s$ and $x(0)$ there must exists a unilateral input $u(t)$ that steers the network towards a state $\bar x$ such that the projection of $\bar x$ on the subspace $\Chi_{\V_s}$ is $\bar x_s$. This is equivalent to the existence of a point $\bar{\bar {x}}\in \C_r(B)$ such that
 \begin{equation}\label{eq:for_proofb}
 \proj_{\Chi_{\V_s}}  \bar{\bar {x}} = \bar x_s - \proj_{\Chi_{\V_s}} \exp(At)x(0).    
 \end{equation}
As from Theorem \ref{thm:2} if $\V_s$ is unilaterally reachable then $\proj_{\Chi_{V_s}} \C_r(B) = {\Chi_{V_s}}$, by hypothesis a point $\bar{\bar{x}}$ fulfilling \eqref{eq:for_proofb} exists for all $\bar x_s$ and $x(0)$. Hence, $\V_s$ is unilaterally controllable.
\end{proof}
The equivalence between unilateral reachability and controllability of node subsets allows to derive a set of corollaries that characterize partial unilateral controllability of network systems.
\begin{mycor}\label{cor:cont}
Given network \eqref{eq:net_eq} and a node subset $\V_s \subseteq \V$, if Assumption \ref{assump:unilateral} holds and $\proj_{\Chi_{\V_s}}(\C_r(B)) = \Chi_{\V_s}$, then $\V_s$ is unilaterally controllable.
\end{mycor}
\begin{proof}
Combining Theorems \ref{thm:2} and \ref{thm:3}, the thesis follows.
\end{proof}
\begin{mycor}\label{cor:uni_net}
Let $\mathcal C_r(B)$ be the
unilaterally reachable set of the pair $(A,B)$. There exists a controllable node subset $\V_s$ such that $|\V_s| \geq |\Chi^l(B)|$.
\label{cor}
\end{mycor}
\begin{proof}
From Theorem \ref{thm:uni_con} we know that if $\C_r(B)$ is the unilaterally reachable set, this implies that $\Chi^l(B)$ is its lineality space of dimension $|\Chi^l(B)|$. Moreover, considering that, for a given subspace $\Chi$ there always exists a set of indices $\mathcal {K}$ and the associated subspace $\bar\Chi=\Span(\cup_{i\in\mathcal K}\{\e_i\})$ such that $\rm{proj}_{\bar \Chi} \Chi = \bar \Chi$, this holds also for $\Chi = \Chi^l(B)$. Therefore, from Theorem 2 there exists a unilaterally reachable node subset $\V_s$ of dimension at least equal to $|\V_s| = |\Chi^l(B)|$ that, from Theorem \eqref{thm:3}, is also unilaterally controllable.
\end{proof}
The next corollary provides a sufficient condition guaranteeing that there exists a unilaterally controllable node subset $\V_s$ that contains a given node $i$. 
\begin{mycor}\label{cor:rays}
Given a node, say $i$, if there exists a quadruplet $\lbrace h, j, k, m \rbrace$ such that $r_{j,k}^\T\e_i>0$, $r_{m,h}^\T\e_i<0$, and $\lbrace r_{j,k}, r_{m,h}\rbrace \in \C_r(B)$, then there exists a unilaterally controllable node subset $\V_s$ such that $i\in \V_s$. 
\begin{proof}
Considering that $\proj_{\e_i}\C_r(B) = \e_i$, the thesis follows from Corollary \ref{cor:cont}.
\end{proof}
\end{mycor}
\begin{myrem}
The mathematical treatment of this section substantially differs from the analyses that are performed when seeking complete unilateral controllability \cite{LinAlt:17}.
First, when the network is not completely unilaterally controllable, the set of reachable states is a convex cone instead of a vector space. Second, we needed to show and consider that, albeit the set of reachable states differs from that of controllable states (Lemma \ref{lem:contr}), the reachable and controllable node subsets do coincide (Theorem \ref{thm:3}). Finally, we had to account for the fact that the number of unilaterally reachable nodes can be larger than the maximal dimension of a  unilaterally reachable subspace, as remarked after Theorem \ref{thm:2}.
\end{myrem}

\section{Greedy Algorithm}\label{sec:greedy}
In what follows, we shall leverage the theoretical findings of Section \ref{sec:results} to design a heuristic for solving problem \eqref{eq:prob_form}. Before illustrating the derivation of the algorithm, we need to introduce the spanning set $\W(B)$ of the lineality space $\Chi^l(B) \in \C_r(B)$ , which can be computed as

\thickmuskip = 0\thickmuskip
\begin{align}\nonumber
	\W(B) =  &\{ r_{i,k}, i,k :  \{r_{i,k},\shortminus r_{i,k}\} \in\, \C_r(B) \land
	\Im(r_{i,k}) = 0 \} \\ &  \cup \,
	\label{eq:spann_lin} 
	\{ \{ \Re(r_{i,k}), \Im(r_{i,k}) \}, i,k : \Im(r_{i,k}) \in \C_r(B) \},
\end{align}
where $i = 1, \dots, \mu, k = 1, \dots, v_i$.
Furthermore, let $\mathcal S$ be the set of all subspaces of $\Real^n$ such that
$\proj_\Chi \mathcal C_r(B)=\Chi, \, \forall \Chi\in\mathcal S.$
From Theorem \ref{thm:uni_con}, solving problem \eqref{eq:prob_form}, that is, finding a maximal unilaterally controllable node subset, is equivalent to finding the matrix $B^\star$ that maximizes the cardinality of the largest subspace in $\mathcal S$. Namely,
\begin{equation}\label{eq:int_opt_prob}
    B^\star:= \argmax_{B}\Big(\max_{\Chi\in\mathcal S} \left|\proj_{\Chi}\C_r(B)\right|\Big).
\end{equation}
Unfortunately, \eqref{eq:int_opt_prob} is a combinatorial problem with time computational complexity of order $O(n!)$
that can only be solved through extensive search, which is unfeasible even for a network of a handful of nodes. 

Since finding an exact solution of \eqref{eq:int_opt_prob} is typically unfeasible, we propose a two-step procedure for the selection of matrix $B$ whose computational complexity is determined by that of finding the Jordan form $J$, that is, $O(n^4)$. In Step 1, 
we seek for a heuristic approach that tries to maximize the lineality $|\Chi^l(B)|$, which from Corollary \ref{cor:uni_net} is a lower bound for the cardinality $|\V_s|$ of the unilaterally controllable node subset $\V_s$. Then, Step 2 attempts to add to $\V_s$ the nodes fulfilling the sufficient condition for node unilateral controllability given in Corollary \ref{cor:rays}.

\subsection*{Step 1. Heuristic maximizing $|\Chi^l(B)|$.}
Here, we seek for the suboptimal solution 
\begin{equation}
    \widetilde B^\star:= \argmax_{B} \ |\Chi^l|, \quad \Chi_l\subseteq \C_r(B)
\end{equation}
to problem \eqref{eq:int_opt_prob}.
The heuristic we propose (Step \ref{alg:step1}) takes as inputs the matrix $A$ and the number of available inputs $m$. Denoting $B_k$ the $\widetilde B$ selected at the $k$-th iteration,
the algorithm starts with $B_0=0_{n\times m}$. Then, at each iteration, one or two columns are added to $\widetilde{B}$.
Defining $\Delta(\beta): = |\Chi^l([B_{k-1}, \beta])|-|\Chi^l(B_{k-1})|$, where $\beta\in\mathcal B$, we can now distinguish two different cases: 
\begin{enumerate}
    \item If there exists $\beta\in\mathcal B$ such that $\Delta(\beta)>0$, a single column is added at step $k$, that is, $ B_k=[B_{k-1},\beta^\star],$ where
    
    \begin{subnumcases}{\beta^\star \shorteq}
        \label{eq:one_input_a}
         \argmax_{\beta\in\mathcal B}|\Chi^l(B_k)|, \text{if } \exists{!}  \argmax_{\beta\in \mathcal{B}}|\Chi^l(B_k)|, 
        \\
        \label{eq:one_input_b}
        \argmax_{\beta \in \mathcal{B}} |\C_{r}(B_k)|,\hspace{0.6 mm}\text{if } \nexists{!}  \argmax_{\beta\in \mathcal{B}}|\Chi^l(B_k)|\footnotemark.
    \end{subnumcases}
    
\addtocounter{footnote}{0}
\footnotetext{With a slight abuse of notation, here we mean that such a $\beta$ exists but is not unique.}
\stepcounter{footnote}
    
 \item If, instead, a $\beta\in\mathcal B$ such that $\Delta(\beta)>0$ does not exist, we add two columns to $B_{k-1}$ at step $k$, that is,
     $B_k=[B_{k-1},\beta^{ \star\star},-\beta^{   \star\star}], $
     where
    \begin{equation}\label{eq:two_inputs}
     \beta^{\star\star}=
     \arg\max_{\beta \in \mathcal{B}} |\Chi^l([B_{k-1},\beta, -\beta])|.
     \end{equation}
\end{enumerate}

Summing up, at each step $k$ our updating rule attempts to add the input that maximizes the lineality $|\Chi^l(B_k)|$. When such an input is not unique, it selects the input that adds the largest number of rays in $\C_r(B_k)$. If instead we cannot find a $\beta$ such that $\Delta(\beta)$ is positive, then we add the two inputs that maximize $|\Chi^l(B_{k})|$. The algorithm stops when $B_k\in\mathbb{R}^{n\times m}$. 
Note that this first step has a computational complexity of $O(n^4)$, due to the evaluation of the Jordan form of $A$. 
\begin{algorithm}[hb!]\caption{Maximizing the lineality $|\Chi^l(B)|$.
 \label{alg:step1}}
	\begin{algorithmic}
	\State \textbf{Inputs:} {$A$, $m$}
	\Procedure{Initialization} {$B_0 = \varnothing$, $\Chi^l(B_0) = \mathcal{O}$,  $\W(B_0) = \varnothing$ 
	}
	\While{$k \leq m -1$} 
		 \If{$\exists \; \beta\in \mathcal{B} : |\Chi^l([B_{k-1}, \; \beta])|>|\Chi^l(B_{k-1})|$} 
    		 \State compute $\beta^\star$ as in \eqref{eq:one_input_a}, \eqref{eq:one_input_b}
    		 \State set $B_k = [B_{k-1}, \ \beta^\star]$
    		 \State compute $\W(B_k)$
    		 \State $k = k + 1$
		 \Else
            \State compute $\beta^{\star\star}$ as in \eqref{eq:two_inputs}
            \State set $B_k = [B_{k-1}, \ \beta^{\star \star}, \ -\beta^{\star \star}]$ and 
            \State compute $\W(B_k)$
            \State $k=k+2$
        \EndIf 
	\EndWhile
	\If{$k = m$}
    		 \State compute $\beta^\star$
    		 as in \eqref{eq:one_input_a}, \eqref{eq:one_input_b}
    		 \State set $B_k=[B_{k-1}, \beta^\star]$
    		 \State compute $\W(B_k)$
    \EndIf
\EndProcedure
\State \textbf{Outputs:} $\widetilde B=B_m$, $\W(\widetilde B)$
	\end{algorithmic}
\end{algorithm}

Once we have computed $\widetilde B = B_m$, we need to identify one of the unilaterally controllable node subsets $\V^1_s$ corresponding to $\widetilde{B}$. To this aim, we leverage Corollary \ref{cor:uni_net}, which states that there exists a unilaterally controllable node subset $\V^1_s$ with $|\V^1_s| =|\Chi^l(\widetilde B)|$ such that $\proj_{\Chi_{\V^1_s}}\Chi^l(\widetilde B) = \Chi_{\V^1_s}$. To identify such a node subset, we compute the set $\W(\widetilde B)$ according to \eqref{eq:spann_lin}. Then, we build the set $\V^1_s$ so that the elements of the sets $\V^1_s$ and $\W(\widetilde B)$ can be associated into $|\V^1_s|$ pairs $(v_j,w_i)$ such that (i) no pairs share a common element and (ii) each pair $(v_j,w_i)$ is such that $\e_{v_j}^{\mathrm{T}}w_i\neq 0$. 
Finding this association can be recast as the problem of finding the maximum matching \cite{edmonds1965maximum} of an unbalanced bipartite graph $\G_b = (\V_b,\E_b)$. Here, $\V_b := \V_w \cup \V_s$ is the set of vertexes and each node in $\V_w$ represents an element of $\W(\widetilde B)$. The set of edges 
$\E_b =\{(i,j)|\e_i^\T w_j\neq 0\land w_j\in \W(\widetilde B)\}$ defines all the possible associations, by appriopriately connecting the nodes in $\V_w$ to those in $\V_s$. Finding a maximum matching is possible by means of the Hopcroft-Karp algorithm \cite{hopcroft1973n} and thus the computational complexity of this sub-step is  $O(\sqrt{|\V_b|}|\E_b|) \leq O(n^{5/2})$.

\subsection*{Step 2. Enlarging the unilaterally controllable node subset.}

In the second step, we enrich the unilaterally controllable node subset $\mathcal V_s^1$ by exploiting the set $\C_r(\widetilde B)\setminus\Chi^l(\widetilde B)$.
To do so, let us define the set
\begin{equation}
\begin{aligned}
    {\label{defQ}}
    \Q(\widetilde B): = 
    &\lbrace r_{i,j}, i,j : r_{i,j}\in \C_r(\widetilde B) \land \shortminus r_{i,j}\notin \C_r(\widetilde B)\rbrace \cup \\
    & \lbrace r_{i,j} , i,j : \shortminus r_{i,j}\in \C_r(\widetilde B) \land r_{i,j} \notin \C_r(\widetilde B)\rbrace,
\end{aligned}
\end{equation}
whose positive span is $\C_r(\widetilde B)\setminus\Chi^l(\widetilde B)$. Then, let us define the matrix $Q$ as the matrix obtained by juxtaposing the elements of $\Q(\widetilde B)$ column-wise. Exploiting Corollary \ref{cor:rays}, we then add a node to $\V_s^1$ whenever the $i$-th row of $Q$ encompasses two nonzero entries, say $q_{ij}$ and $q_{im}$, that are such that $\sgn(q_{ij})=-\sgn(q_{im})$. Let us note that the computational complexity of this step of the algorithm is $O(n^3)$.

\begin{algorithm}[hb!]\caption{Enlarging the unilaterally controllable node subset associated to $\widetilde B$
 \label{alg:step3}}
	\begin{algorithmic}
	\State \textbf{Inputs:} {$Q,\V^1_s$}
	\For{i = 1, \dots, n}
        \If {$\exists \ l,m : \sgn(q_{i,m}) = -\sgn(q_{i,l})$ }
        \State Set $\V_s^1 = \V_s^1 \cup \{i\}$
        \State Remove the $l$-th and $m$-th columns from $Q$
        \EndIf
        \EndFor
    \State \textbf{Output:} $\V_s = \V_s^1.$
	\end{algorithmic}
\end{algorithm}

\subsection*{Application on a sample network}
To illustrate our heuristic, we consider a linear network dynamical system on the graph $\mathcal{G}$ depicted in Fig. 1, whose dynamics is described by matrix
\[
A =
\begin{bmatrix}
1	&-4	    &0	&0	&\ \ 0	&\ \ 0	&\ \ 0 \\
4	&\ \ 1	&0  &0	&\ \ 0	&\ \ 0	&\ \ 0 \\
1	&\ \ 0	&3	&0	   &-1 	&\ \ 0	   &-1 \\
0	&\ \ 0	&1	&4	&\ \ 1	&\ \ 0	&\ \ 4 \\
0	&\ \ 0	&0	&0	&\ \ 2	   &-3	&\ \ 0 \\
0	&\ \ 0	&0	&0	&\ \ 3	&\ \ 2	&\ \ 0 \\
0	&\ \ 0	&0	&0	   &-3	&\ \ 0  &\ \ 0 \\
\end{bmatrix},
\]
with spectrum $\{ 4, 3, 0, 1+4i, 1-4i, 2+3i, 2-3i \}.$
Let us assume that we can inject $m=2$ unilateral controls.
The input matrix $\widetilde B$ is designed following Step \ref{alg:step1}, that is, by maximizing the lineality $|\Chi^l(\widetilde B)|$. At time instant $k=1$, four possible selections of $\beta$ ($\e_2$, $\e_6$, $-\e_2$ and $-\e_6$) yield the same (positive) $\Delta(\beta)$. Hence, $\beta^\star$ should be selected among these four according to \eqref{eq:one_input_b}. However, since all choices would yield the same $|\mathcal{C}(B_1)|$, the selection is performed randomly, and we elect $B_1 = -\e_6$, with the set $\W(B_1)$ being $[r_6, r_7]$. 
At $k = 2$, $-\e_2$ is the unique $\beta$ returning $\Delta(\beta)>0$. Hence, we select node $2$ as the second and last node where a control signal is injected, i.e., we set $\widetilde B = B_2 = [-\e_6; -\e_2]$ and $\W(\widetilde{B}) = [r_1, r_4, r_5, r_6, r_7]$. 

Having selected the matrix $B$, we now turn to finding one of the possibly multiple unilaterally controllable node subsets $\V^1_s$ such that $|\V^1_s| = |\Chi^l(\widetilde B)|$ by solving the maximum matching problem.
Among the multiple equivalent solutions to this problem, we randomly pick $\V^1_s=\{v_1, v_2, v_4, v_5, v_6\}$. Finally, we compute, from \eqref{defQ}, $Q(\widetilde B)=\{-r_2,r_3\}$, and from Step \ref{alg:step3} of the proposed heuristic we can enlarge the unilaterally controllable node subset with node $3$, that is $\V_s = \{v_1,v_2, v_3, v_4, v_5, v_6\}$.  Interestingly, in this simple example, we find that $|\V_s|>|\Chi^l|$, that is, the number of unilaterally controllable nodes is greater than the largest unilaterally controllable subspace.
\begin{figure}[ht!]
    \centering \vspace{0.4cm}
\includegraphics[scale = 0.6]{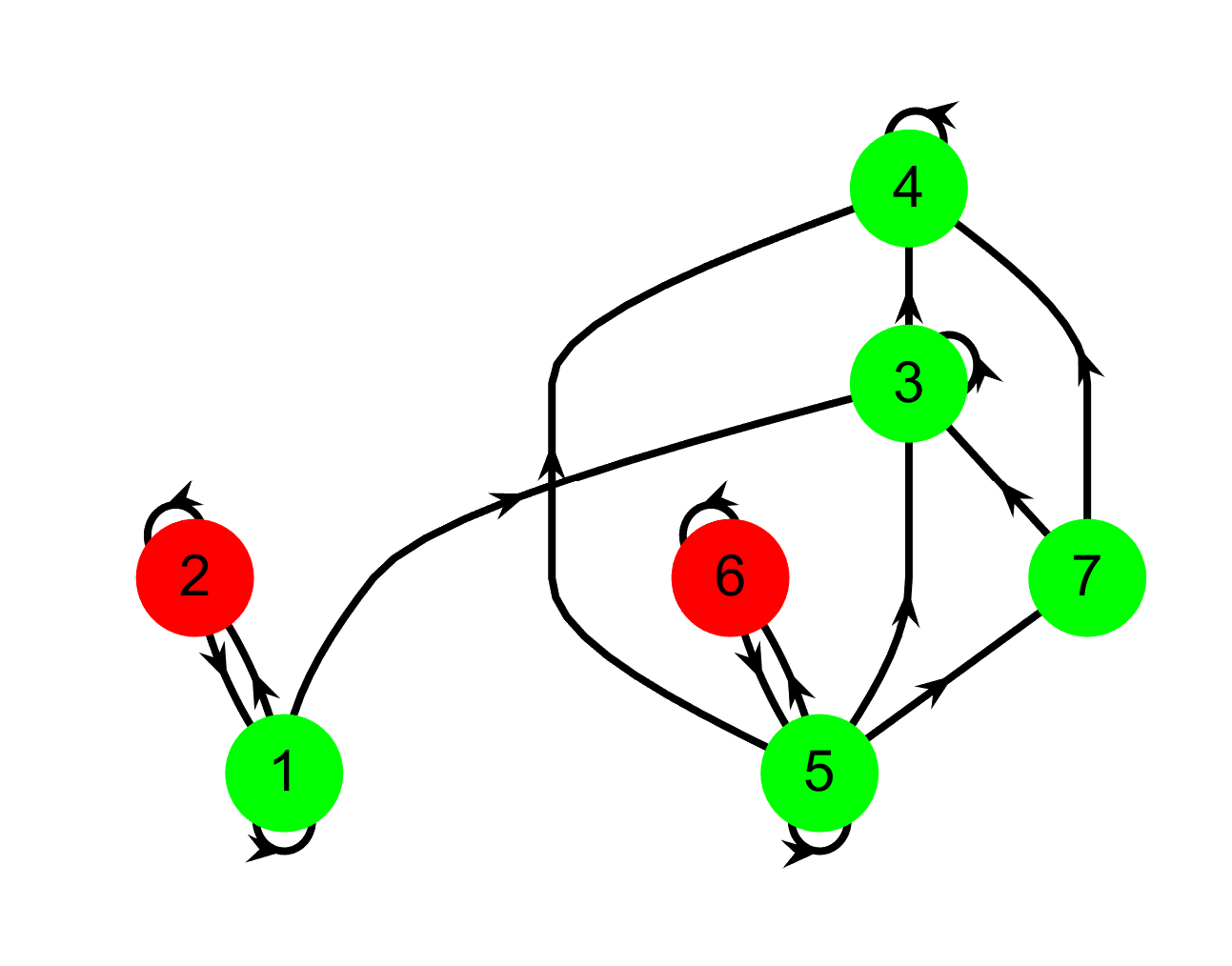}
 \caption{Network topology: the nodes in red are the driver nodes where, according to the proposed heuristic, the negative unilateral control inputs are injected.}
    \label{fig:mesh1}
\end{figure}
\section{Conclusions}
In this letter, we have studied controllability of linear network dynamical systems when the inputs are unilateral. Specifically, we focused on the case where the constraint on the number and type of inputs prevents the achievement of complete controllability of the network system, whereby only a node subset can be made controllable. In this setting, we have identified conditions for unilateral reachability and  controllability of a node set, which we found to be equivalent, different from the general case of subsystems, where we have proved that reachability does not imply controllability. After showing that  maximizing the size of a controllable nodes subset is a combinatorial problem, we have leveraged the theoretical findings on unilateral controllability to build an heuristic that can find a suboptimal solution to this problem in polynomial time, as illustrated on a sample network. 

Our work has laid the foundations of partial controllability under unilateral inputs, thus paving the ways for future studies in this area of research. First, alternative heuristic approaches may be developed and tested against the one proposed in this manuscript. Moreover, once partial controllability has been guaranteed, the problem of evaluating the energy associated to the control action arises, thereby minimum energy control problems could be formulated in this setting.

\section*{Acknowledgments}
The authors wish to thank Professor Claudio Altafini of the Dept. of Electrical Engineering,
Linkoping University, for taking part to insightful discussions on this topic.

\bibliographystyle{IEEEtran}
\bibliography{ref}
\end{document}